\documentclass{amsart}
\usepackage{amsmath,amscd,amsthm,amsfonts}

\theoremstyle{plain}
\newtheorem{theorem}                 {Theorem}      [section]

\newtheorem{lemma}        [theorem]  {Lemma}

\theoremstyle{definition}
\newtheorem{example}      [theorem]  {Example}

\newtheorem{definition}   [theorem]  {Definition}

\numberwithin{equation}{section}

\begin{document}
\baselineskip 18pt \larger

\def \theo-intro#1#2 {\vskip .25cm\noindent{\bf Theorem #1\ }{\it #2}}

\newcommand{\trace}{\operatorname{trace}}

\def \dn{\mathbb D}
\def \nn{\mathbb N}
\def \zn{\mathbb Z}
\def \qn{\mathbb Q}
\def \rn{\mathbb R}
\def \cn{\mathbb C}
\def \hn{\mathbb H}
\def \P{\mathbb P}
\def \can{Ca}

\def \S{\mathcal S}
\def \A{\mathcal A}
\def \B{\mathcal B}
\def \C{\mathcal C}
\def \E{\mathcal E}
\def \F{\mathcal F}
\def \G{\mathcal G}
\def \H{\mathcal H}
\def \I{\mathcal I}
\def \L{\mathcal L}
\def \M{\mathcal M}
\def \N{\mathcal N}
\def \O{\mathcal O}
\def \P{\mathcal P}
\def \Q{\mathcal Q}
\def \R{\mathcal R}
\def \V{\mathcal V}
\def \W{\mathcal W}

\def\Re{\mathfrak R\mathfrak e}
\def\Im{\mathfrak I\mathfrak m}
\def\Co{\mathfrak C\mathfrak o}
\def\Or{\mathfrak O\mathfrak r}

\def \ip #1#2{\langle #1,#2 \rangle}
\def \spl#1#2{( #1,#2 )}

\def \lb#1#2{[#1,#2]}

\def \pror#1{\rn P^{#1}}
\def \proc#1{\cn P^{#1}}
\def \proh#1{\hn P^{#1}}

\def \gras#1#2{G_{#1}(\cn^{#2})}

\def \b{\mathfrak{b}}
\def \g{\mathfrak{g}}
\def \h{\mathfrak{h}}
\def \k{\mathfrak{k}}
\def \m{\mathfrak{m}}
\def \p{\mathfrak{p}}
\def \q{\mathfrak{q}}
\def \r{\mathfrak{r}}
\def \un{\mathfrak{u}}

\def \GLR#1{\text{\bf GL}_{#1}(\rn)}
\def \GLRP#1{\text{\bf GL}^+_{#1}(\rn)}
\def \glr#1{\mathfrak{gl}_{#1}(\rn)}
\def \GLC#1{\text{\bf GL}_{#1}(\cn)}
\def \glc#1{\mathfrak{gl}_{#1}(\cn)}
\def \GLH#1{\text{\bf GL}_{#1}(\hn)}
\def \glh#1{\mathfrak{gl}_{#1}(\hn)}
\def \GLD#1{\text{\bf GL}_{#1}(\dn)}
\def \gld#1{\mathfrak{gl}_{#1}(\dn)}

\def \SLR#1{\text{\bf SL}_{#1}(\rn)}
\def \slr#1{\mathfrak{sl}_{#1}(\rn)}
\def \SLC#1{\text{\bf SL}_{#1}(\cn)}
\def \slc#1{\mathfrak{sl}_{#1}(\cn)}

\def \O#1{\text{\bf O}(#1)}
\def \SO#1{\text{\bf SO}(#1)}
\def \so#1{\mathfrak{so}(#1)}
\def \SOs#1{\text{\bf SO}^*(#1)}
\def \sos#1{\mathfrak{so}^*(#1)}
\def \SOO#1#2{\text{\bf SO}(#1,#2)}
\def \soo#1#2{\mathfrak{so}(#1,#2)}
\def \SOC#1{\text{\bf SO}(#1,\cn)}
\def \SOc#1{\text{\bf SO}(#1,\cn)}
\def \soc#1{\mathfrak{so}(#1,\cn)}

\def \SUs#1{\text{\bf SU}^*(#1)}
\def \sus#1{\mathfrak{su}^*(#1)}
\def \Us#1{\text{\bf U}^*(#1)}
\def \us#1{\mathfrak{u}^*(#1)}

\def \U#1{\text{\bf U}(#1)}
\def \u#1{\mathfrak{u}(#1)}
\def \UU#1#2{\text{\bf U}(#1,#2)}
\def \uu#1#2{\mathfrak{u}(#1,#2)}
\def \SU#1{\text{\bf SU}(#1)}
\def \su#1{\mathfrak{su}(#1)}
\def \SUU#1#2{\text{\bf SU}(#1,#2)}
\def \suu#1#2{\mathfrak{su}(#1,#2)}

\def \Sp#1{\text{\bf Sp}(#1)}
\def \sp#1{\mathfrak{sp}(#1)}
\def \Spp#1#2{\text{\bf Sp}(#1,#2)}
\def \spp#1#2{\mathfrak{sp}(#1,#2)}
\def \SpR#1{\text{\bf Sp}(#1,\rn)}
\def \spR#1{\mathfrak{sp}(#1,\rn)}
\def \SpC#1{\text{\bf Sp}(#1,\cn)}
\def \spc#1{\mathfrak{sp}(#1,\cn)}

\def \d#1{\mathfrak{d}(#1)}
\def \s#1{\mathfrak{s}(#1)}
\def \sym#1{\text{Sym}(\rn^{#1})}
\def \symc#1{\text{Sym}(\cn^{#1})}

\def \gradh#1{\text{grad}_{\H}(#1 )}
\def \gradv#1{\text{grad}_{\V}(#1 )}

\def \nab#1#2{\hbox{$\nabla$\kern -.3em\lower 1.0 ex
    \hbox{$#1$}\kern -.1 em {$#2$}}}

\allowdisplaybreaks
\title{Harmonic morphisms from the classical
\\non-compact semisimple Lie groups\\(version 2.043)}

\author{Sigmundur Gudmundsson}
\author{Anna Sakovich}



\address
{Mathematics, Faculty of Science, Lund University, Box 118, S-221
00 Lund, Sweden} \email{Sigmundur.Gudmundsson@math.lu.se}

\address
{Faculty of Pre-University Education, Belarusian State
University, Oktyabrskaya Str. 4, Minsk 220030, Belarus}
\email{anya\_sakovich@tut.by}

\begin{abstract}
We construct the first known complex valued harmonic morphisms
from the non-compact Lie groups $\SLR n$, $\SUs{2n}$ and $\SpR n$
equipped with their standard Riemannian metrics.  We then
introduce the notion of a {\it bi-eigenfamily} and employ this to
construct the first known solutions on the non-compact Riemannian
$\SOs{2n}$, $\SOO pq$, $\SUU pq$ and $\Spp pq$. Applying a {\it
duality principle} we then show how to manufacture the first known
complex valued harmonic morphisms from the compact Lie groups $\SO
n$, $\SU n$ and $\Sp n$ equipped with semi-Riemannian metrics.
\end{abstract}

\maketitle

\section{Introduction}

In differential geometry the notion of a minimal submanifold of a
given ambient space is of great importance. Harmonic morphisms
$\phi:(M,g)\to(N,h)$ between semi-Riemannian manifolds are useful
tools for the construction of such objects. They are solutions to
over-determined non-linear systems of partial differential
equations determined by the geometric data of the manifolds
involved. For this reason harmonic morphisms are difficult to find
and have no general existence theory, not even locally.

For the existence of harmonic morphisms $\phi:(M,g)\to (N,h)$ it
is an advantage that the target manifold $N$ is a surface i.e. of
dimension $2$. In this case the problem is invariant under
conformal changes of the metric on $N^2$. Therefore, at least for
local studies, the codomain can be taken to be the complex plane
with its standard flat metric.   For the general theory of
harmonic morphisms between semi-Riemannian manifolds we refer to
the excellent book \cite{Bai-Woo-book} and the regularly updated
on-line bibliography \cite{Gud-bib}.

In \cite{Gud-Sak-1} the current authors introduce the notion of an
{\it eigenfamily} of complex valued functions from semi-Riemannian
manifolds. This is used to manufacture a variety of locally
defined harmonic morphisms from the classical compact Lie groups
$$\SO n,\ \ \SU n\ \ \text{and}\ \ \Sp n$$ equipped with their standard
Riemannian metrics. A general {\it duality principle} is developed
and employed to construct solutions on the non-compact Lie groups

$$\SLR n,\ \ \SUs{2n},\ \ \SpR n,$$

$$\SOs{2n},\ \ \SOO pq,\ \ \SUU pq\ \ \text{and}\ \ \Spp pq$$

\noindent equipped with their standard dual semi-Riemannian
metrics.

The current paper is devoted to the study of the {\it dual
problem}. We are now mainly interested in the classical
non-compact semisimple Lie groups equipped with their standard
Riemannian metrics. We construct the first known locally defined
complex valued harmonic morphisms from
$$\SLR n,\ \ \SUs{2n}\ \ \text{and}\ \ \SpR n$$
which are not invariant under the action of the subgroups $\SO n$,
$\Sp n$ or $\U n$, respectively, see \cite{Gud-Sve-3}. We
generalize the idea of an eigenfamily and introduce the notion of
a {\it bi-eigenfamily} of complex valued functions on a
semi-Riemannian manifold. This leads to the construction of the
first known locally defined complex valued harmonic morphisms on
the important non-compact Riemannian Lie groups $$\SOs{2n},\ \
\SOO pq,\ \ \SUU pq\ \ \text{and}\ \ \Spp pq.$$ Employing the
earlier mentioned {\it duality principle} we then show how to
produce the first known locally defined complex valued harmonic
morphisms on the compact Lie groups $$\SO n,\ \ \SU n\ \
\text{and}\ \ \Sp n$$ equipped with their standard dual
semi-Riemannian metrics. It should be noted that the non-compact
semisimple Lie groups
$$\SOC n,\ \ \SLC n\ \ \text{and}\ \ \SpC n$$ are complex manifolds and
hence their coordinate functions form orthogonal harmonic
families, see \cite{Gud-Sak-1}.  This means that in these cases
the problem is more or less trivial.

Throughout this article we assume, when not stating otherwise,
that all our objects such as manifolds,  maps etc. are smooth i.e.
in the $C^{\infty}$-category. For our notation concerning Lie
groups we refer to the wonderful book \cite{Kna}.

\section{Harmonic Morphisms}

Let $M$ and $N$ be two manifolds of dimensions $m$ and $n$,
respectively. Then a semi-Riemannian metric $g$ on $M$ gives rise
to the notion of a Laplacian on $(M,g)$ and real-valued harmonic
functions $f:(M,g)\to\rn$. This can be generalized to the concept
of a harmonic map $\phi:(M,g)\to (N,h)$ between semi-Riemannian
manifolds being a solution to a semi-linear system of partial
differential equations, see \cite{Bai-Woo-book}.

\begin{definition}
A map $\phi:(M,g)\to (N,h)$ between semi-Riemannian manifolds is
called a {\it harmonic morphism} if, for any harmonic function
$f:U\to\rn$ defined on an open subset $U$ of $N$ with
$\phi^{-1}(U)$ non-empty, the composition
$f\circ\phi:\phi^{-1}(U)\to\rn$ is a harmonic function.
\end{definition}

The following characterization of harmonic morphisms between
semi-Riemannian manifolds is due to Fuglede and generalizes the
corresponding well-known result of \cite{Fug-1,Ish} in the
Riemannian case.  For the definition of horizontal conformality we
refer to \cite{Bai-Woo-book}.

\begin{theorem}\cite{Fug-2}
  A map $\phi:(M,g)\to (N,h)$ between semi-Rie\-mannian manifolds is a
  harmonic morphism if and only if it is a horizontally (weakly)
  conformal harmonic map.
\end{theorem}

The next result generalizes the corresponding well-known theorem
of Baird and Eells in the Riemannian case, see \cite{Bai-Eel}. It
gives the theory of harmonic morphisms a strong geometric flavour
and shows that the case when $n=2$ is particularly interesting. In
that case the conditions characterizing harmonic morphisms are
independent of conformal changes of the metric on the surface
$N^2$.  For the definition of horizontal homothety we refer to
\cite{Bai-Woo-book}.

\begin{theorem}\cite{Gud-1}\label{theo:semi-B-E}
Let $\phi:(M,g)\to (N^n,h)$ be a horizontally conformal submersion
from a semi-Riemannian manifold $(M,g)$ to a Riemannian manifold
$(N,h)$. If
\begin{enumerate}
\item[i.] $n=2$ then $\phi$ is harmonic if and only if $\phi$ has
minimal fibres, \item[ii.] $n\ge 3$ then two of the following
conditions imply the other:
\begin{enumerate}
\item $\phi$ is a harmonic map, \item $\phi$ has minimal fibres,
\item $\phi$ is horizontally homothetic.
\end{enumerate}
\end{enumerate}
\end{theorem}

In what follows we are mainly interested in complex valued
functions $\phi,\psi:(M,g)\to\cn$ from semi-Riemannian manifolds.
In this situation the metric $g$ induces the complex-valued
Laplacian $\tau(\phi)$ and the gradient $\text{grad}(\phi)$ with
values in the complexified tangent bundle $T^{\cn}M$ of $M$.  We
extend the metric $g$ to be complex bilinear on $T^{\cn} M$ and
define the symmetric bilinear operator $\kappa$ by
$$\kappa(\phi,\psi)= g(\text{grad}(\phi),\text{grad}(\psi)).$$ Two
maps $\phi,\psi: M\to\cn$ are said to be {\it orthogonal} if
$$\kappa(\phi,\psi)=0.$$  The harmonicity and horizontal
conformality of $\phi:(M,g)\to\cn$ are given by the following
relations
$$\tau(\phi)=0\ \ \text{and}\ \ \kappa(\phi,\phi)=0.$$

\begin{definition}\label{defi:eigen}
Let $(M,g)$ be a semi-Riemannian manifold.  Then a set
$$\E =\{\phi_i:M\to\cn\ |\ i\in I\}$$ of complex valued functions
is said to be an {\it eigenfamily} on $M$ if there exist complex
numbers $\lambda,\mu\in\cn$ such that
$$\tau(\phi)=\lambda\phi\ \ \text{and}\ \ \kappa
(\phi,\psi)=\mu\phi\psi$$ for all $\phi,\psi\in\E$.
\end{definition}

The next result shows that an eigenfamily on a semi-Riemannian
manifold can be used to produce a variety of local harmonic
morphisms.

\begin{theorem}\cite{Gud-Sak-1}\label{theo:rational}
Let $(M,g)$ be a semi-Riemannian manifold and $$\E =\{\phi_1,\dots
,\phi_n\}$$ be a finite eigenfamily of complex valued functions on
$M$. If $P,Q:\cn^n\to\cn$ are linearily independent homogeneous
polynomials of the same positive degree then the quotient
$$\frac{P(\phi_1,\dots ,\phi_n)}{Q(\phi_1,\dots ,\phi_n)}$$ is a
non-constant harmonic morphism on the open and dense subset
$$\{p\in M| \ Q(\phi_1(p),\dots ,\phi_n(p))\neq 0\}.$$
\end{theorem}

\section{The Riemannian Lie group $\GLC n$}

Let $G$ be a Lie group with Lie algebra $\g$ of left-invariant
vector fields on $G$.  Then a Euclidean scalar product $g$ on the
algebra $\g$ induces a left-invariant Riemannian metric on the
group $G$ and turns it into a Riemannian manifold. If $Z$ is a
left-invariant vector field on $G$ and $\phi,\psi:U\to\cn$ are two
complex valued functions defined locally on $G$ then the first and
second order derivatives satisfy
$$Z(\phi)(p)=\frac {d}{ds}[\phi(p\cdot\exp(sZ))]\big|_{s=0},$$
$$Z^2(\phi)(p)=\frac {d^2}{ds^2}[\phi(p\cdot\exp(sZ))]\big|_{s=0}.$$
The tension field $\tau(\phi)$ and the $\kappa$-operator
$\kappa(\phi,\psi)$ are then given by
$$\tau(\phi)=\sum_{Z\in\B}Z^2(\phi),\ \
\kappa(\phi,\psi)=\sum_{Z\in\B}Z(\phi)Z(\psi)$$ where $\B$ is any
orthonormal basis of the Lie algebra $\g$.

Let $\GLC n$ be the complex general linear group equipped with its
standard Riemannian metric induced by the Euclidean scalar product
on the Lie algebra $\glc n$ given by
$$g(Z,W)=\Re\trace ZW^*.$$  For $1\le i,j\le n$ we shall by
$E_{ij}$ denote the element of $\glr n$ satisfying
$$(E_{ij})_{kl}=\delta_{ik}\delta_{jl}$$ and by $D_t$ the diagonal
matrices $$D_t=E_{tt}.$$ For $1\le r<s\le n$ let $X_{rs}$ and
$Y_{rs}$ be the matrices satisfying
$$X_{rs}=\frac 1{\sqrt 2}(E_{rs}+E_{sr}),\ \ Y_{rs}=\frac
1{\sqrt 2}(E_{rs}-E_{sr}).$$  With the above notation we have the
following easily verified matrix identities
$$\sum_{r<s}X_{rs}^2=\frac {(n-1)}2I_n,\ \ \
\sum_{r<s}Y_{rs}^2=-\frac {(n-1)}2I_n,\ \ \
\sum_{t=1}^nD_t^2=I_n,$$

$$\sum_{r<s}X_{rs}E_{jl}X^t_{rs}=\frac 12(E_{lj}
+\delta_{lj}(I_n-2E_{lj})),$$

$$\sum_{r<s}Y_{rs}E_{jl}Y^t_{rs}=-\frac 12( E _{lj}
-\delta _{lj}I_n),$$

$$\sum_{t=1}^nD_{t}E_{jl}D^t_{t}=\
\delta _{jl} E _{lj}.$$

\section{The Riemannian Lie group $\SLR n$}

In this section we construct eigenfamilies of complex valued
functions on the general linear group $\GLR n$.  These families
can be used to construct local harmonic morphisms on the
non-compact semisimple special linear group
$$\SLR n=\{ x\in\GLR n|\ \det x = 1\}.$$

The Lie group $\rn^*$ of non-zero real numbers acts on $\GLR n$ by
multiplication $(r,x)\mapsto rx$ and the orbit space of this
action is $\SLR n$. The Lie algebra $\glr n$ of $\GLR n$ consists
of all real $n\times n$ matrices and for this we have the
canonical orthonormal basis
$$\{X_{rs}, Y_{rs}|\ 1\le r<s\le n\}\cup\{D_t|\ t=1,\dots ,n\}.$$

\begin{lemma}\label{lemm:special}
For $1\le i,j\le n$ let $x_{ij}:\GLR n\to\rn$ be the real valued
coordinate functions given by
$$x_{ij}:x\mapsto e_i\cdot x\cdot e_j^t$$ where $\{e_1,\dots
,e_n\}$ is the canonical basis for $\rn^n$. Then the following
relations hold  $$\tau(x_{ij})=x_{ij},$$
$$\kappa(x_{ij},x_{kl})=\delta_{jl}\sum_{t=1}^nx_{it}x_{kt}.$$
\end{lemma}

\begin{proof}  It follows directly from the definition of the functions
$x_{ij}$ that if $X$ is an element of the Lie algebra $\glr n$
then the first and second order derivatives satisfy
$$X(x_{ij}):x\mapsto e_i\cdot x\cdot X\cdot e_j^t\ \ \text{and}
\ \  X^2(x_{ij}):x\mapsto e_i\cdot x\cdot X^2\cdot e_j^t.$$
Employing the above mentioned matrix identities we then yield

\begin{eqnarray*}
\tau(x_{ij})&=&\sum_{r<s}X_{rs}^2(x_{ij})+\sum_{r<s}Y_{rs}^2(x_{ij})
+\sum_{t=1}^nD_t^2(x_{ij})\\
&=&e_i\cdot x\cdot
(\sum_{r<s}X_{rs}^2+\sum_{r<s}Y_{rs}^2+\sum_{t=1}^nD_t^2)\cdot e_j^t\\
&=&x_{ij}.
\end{eqnarray*}

\begin{eqnarray*}
\kappa(x_{ij},x_{kl})
&=&\sum_{r<s}X_{rs}(x_{ij})X_{rs}(x_{kl})+\sum_{r<s}Y_{rs}(x_{ij})Y_{rs}(x_{kl})\\
& & \quad +\sum_{t=1}^nD_t(x_{ij})D_t(x_{kl})\\
&=&e_i\cdot x\cdot\big(\sum_{r<s}X_{rs}\cdot E_{jl}\cdot
X_{rs}^t\big)\cdot x^t\cdot e_k^t\\
& &\quad +e_i\cdot x\cdot\big(\sum_{r<s}Y_{rs}\cdot E_{jl}\cdot
Y_{rs}^t\big)\cdot x^t\cdot e_k^t\\
& &\quad +e_i\cdot x\cdot\big(\sum_{t=1}D_t\cdot E_{jl}\cdot
D_t^t\big)\cdot x^t\cdot e_k^t\\
&=&\delta_{jl}\sum_{t=1}^nx_{it}x_{kt}.
\end{eqnarray*}
\end{proof}

Let $P,Q:\GLR n\to\cn$ be homogeneous polynomials of the
coordinate functions $x_{ij}:\GLR n\to\cn$ of degree one i.e. of
the form
$$P=\trace (Ax^t)=\sum_{i,j=1}^n a_{ij}x_{ij}\ \ \text{and}
\ \ Q=\trace (Bx^t)=\sum_{k,l=1}^n b_{kl}x_{kl}$$ for some
$A,B\in\cn^{n\times n}$.  Then it is easily seen that
$$\kappa (P,Q)=\trace (x^tAB^tx)=\trace (AB^txx^t).$$

\begin{theorem}\label{theo:special}  Let $V$ be a maximal
isotropic subspace of $\cn^n$ and $M(V)$ be the set of matrices
with rows all contained in $V$. Then the complex vector space
$$\E_V=\{\phi _A:\GLR n\to\cn\ |\ \phi_A(x)=
\trace (Ax^t),\ A\in M(V)\} $$ is an eigenfamily on $\GLR n$.
\end{theorem}

It should be noted that the local harmonic morphisms on $\GLR n$
obtained by Theorem \ref{theo:rational} and Theorem
\ref{theo:special} are invariant under the action of $\rn^*$ and
hence induce local harmonic morphisms on the special linear group
$\SLR n$.

\begin{example}
The $3$-dimensional Lie group $\SLR 2$ is given by
$$\SLR 2=\{\begin{pmatrix} a& b
\\ c & d\end{pmatrix}\in\rn^{2\times 2} |\ ad-bc=1\}.$$
The vector $(1,i)$ generates a $1$-dimensional maximal isotropic
subspace $V$ of $\cn^2$ inducing the eigenfamily
$$\E_V=\{\phi_A:\GLR 2\to\cn\ |\ A\in M(V)\}$$
on $\GLR 2$.  By choosing $$A=\begin{pmatrix} 1& i
\\ 0 & 0\end{pmatrix}, \ \ B=\begin{pmatrix} 0& 0
\\ 1 & i\end{pmatrix}$$
and applying Theorem \ref{theo:rational} we yield the well known
globally defined harmonic morphism $\phi={\phi_A}/{\phi_B}:\SLR
2\to H^2$ with
$$\phi(\begin{pmatrix} a& b \\c & d\end{pmatrix})=\frac {a+ib}{c+id}.$$
Here $H^2\cong\SLR 2/\SO 2$ is the hyperbolic upper half plane in
$\cn$.
\end{example}

\section{The Riemannian Lie group $\SUs{2n}$}

In this section we construct eigenfamilies on the non-compact Lie
group
$$\Us{2n}=\{\begin{pmatrix} z & w
\\ -\bar w & \bar z\end{pmatrix}
\ |\ z,w\in\GLC {n}\}.$$  These families can be used to construct
local harmonic morphisms on the semisimple $\SUs{2n}$. The Lie
group $\rn^+$ of positive real numbers acts on $\Us{2n}$ by
multiplication $(r,q)\mapsto rq$ and the orbit space of this
action is
$$\SUs{2n}=\Us{2n}\cap\SLC {2n}.$$
The Lie algebra $\us{2n}$ of $\Us{2n}$ is given by
$$\us{2n}=\{\begin{pmatrix}Z & W \\
-\bar W & \bar Z\end{pmatrix}|\ Z,W \in \glc {n}\}$$ and for this
we have the standard orthonormal basis consisting of the following
elements
$$\frac 1{\sqrt 2}\begin{pmatrix}X_{rs} & 0 \\
0 & X_{rs}\end{pmatrix},\frac 1{\sqrt 2}\begin{pmatrix}Y_{rs} & 0 \\
0 & Y_{rs}\end{pmatrix},\frac 1{\sqrt 2}\begin{pmatrix}D_t & 0 \\
0 & D_t\end{pmatrix},$$

$$\frac 1{\sqrt 2}\begin{pmatrix}iX_{rs} & 0 \\
0 & -iX_{rs}\end{pmatrix},\frac 1{\sqrt 2}\begin{pmatrix}iY_{rs} & 0 \\
0 & -iY_{rs}\end{pmatrix},\frac 1{\sqrt 2}\begin{pmatrix}i D_t & 0 \\
0 & -i D_t\end{pmatrix},$$

$$\frac 1{\sqrt 2}\begin{pmatrix}0 & X_{rs} \\
-X_{rs} & 0\end{pmatrix},\frac 1{\sqrt 2}\begin{pmatrix}0 & Y_{rs} \\
-Y_{rs} & 0\end{pmatrix},\frac 1{\sqrt 2}\begin{pmatrix}0 & D_t \\
-D_t & 0\end{pmatrix},$$

$$\frac 1{\sqrt 2}\begin{pmatrix}0 & iX_{rs} \\
iX_{rs} & 0\end{pmatrix},\frac 1{\sqrt 2}\begin{pmatrix}0 & iY_{rs} \\
iY_{rs} & 0\end{pmatrix},\frac 1{\sqrt 2}\begin{pmatrix}0 & i D_t \\
i D_t & 0\end{pmatrix}$$ where $1\leq r<s \leq n$ and $1\le t\le
n$.

\begin{lemma}\label{lemm:U*(2n)}
For $1\le i,j\le n$ let $z_{ij},w_{ij}:\Us{2n}\to\cn$ be the
complex valued coordinate functions given by
$$z_{ij}:q\mapsto e_i\cdot q\cdot e_j^t,\ \
w_{ij}:q\mapsto e_i\cdot q\cdot e_{n+j}^t$$ where $\{e_1,\dots
,e_{2n}\}$ is the canonical orthonormal basis for $\cn^{2n}$. Then
the following relations hold
$$\tau(z_{ij})=-z_{ij},\ \ \tau(w_{ij})=-w_{ij},$$
$$\kappa(z_{ij}, z_{kl})= 0, \ \ \kappa(w_{ij}, w_{kl})= 0,$$
$$\kappa(z_{ij}, w_{kl})= \delta _{jl}\cdot\sum _{t=1}^n( z_{it}
w_{kt}-w_{it} z_{kt}).$$
\end{lemma}
\begin{proof}  The method of proof is exactly the same as for
Lemma \ref{lemm:special}.
\end{proof}

Let $P,Q:\Us{2n}\to\cn$ be homogeneous polynomials of the
coordinate functions $z_{ij},w_{ij}:\GLR n\to\cn$ of degree one
i.e. of the form
$$P(q)=\trace (Az^t+Bw^t),\ \ Q(q)=\trace (Cz^t+Dw^t)$$
for some $A,B,C,D\in\cn^{n\times n}$. Employing the results of
Lemma \ref{lemm:U*(2n)} a simple calculation shows that
$$\kappa (P,Q)=\trace (((AD^t-BC^t)-(AD^t-BC^t)^t)zw^t).$$
Comparing coefficients we see that
$\kappa(P,P)=\kappa(P,Q)=\kappa(Q,Q)=0$ if and only if the
matrices $$AB^t,\ \ (AD^t-BC^t),\ \ CD^t$$ are symmetric.

\begin{theorem}\label{theo:U*(2n)}
Let $M$ be a non-empty subset of $\cn^{n\times
n}\times\cn^{n\times n}$ such that if $(A,B),(C,D)\in M$ then the
matrices $AB^t,(AD^t-BC^t),CD^t$ are symmetric. Then the complex
vector space
$$\E_M=\{\phi:\Us{2n}\to\cn\ |\ \phi(q)=
\trace (Az^t + Bw^t),\ (A,B)\in M\} $$ is an eigenfamily on
$\Us{2n}$.
\end{theorem}

It is easy to construct non-empty subsets $M$ of $\cn^{n\times
n}\times\cn^{n\times n}$ satisfying the conditions of Theorem
\ref{theo:U*(2n)}.  Here we present two examples.

\begin{example} If $\xi$ is a complex number, then the complex
$n^2$-dimensional vector space
$$M_\xi=\{(A,\xi A)\in\cn^{n\times n}\times\cn^{n\times n}|\
A\in\cn^{n\times n}\}$$ satisfies the conditions above.
\end{example}

\begin{example} If $p$ is an element of $\cn^n$,
then the $2n$-dimensional vector space
$$M_p=\{(pa^t,pb^t)\in\cn^{n\times n}\times\cn^{n\times n}|\
a,b\in\cn^n\}$$ satisfies the conditions above.
\end{example}

The local harmonic morphisms on U*(2n) obtained by Theorem
\ref{theo:rational} and Theorem \ref{theo:U*(2n)} are invariant
under the action of $\rn^+$ and hence induce local harmonic
morphisms on the semisimple Lie group $\SUs{2n}$.

\section{The Riemannian Lie group $\SpR n$}

In this section we construct eigenfamilies on the non-compact
semisimple Lie group
$$\SpR n =\{g\in\SLR{2n}\ |\ g\cdot J_{n}\cdot
g^t=J_{n}\}.$$ The Lie algebra $\spR{n}$ of $\SpR n$ is given by
$$\spR{n}=\{\begin{pmatrix}X & Y \\
Z & -X^t\end{pmatrix}|\ X,Y,Z \in\rn^{n\times n},\ Y=Y^t\
\text{and}\ Z=Z^t\}$$ and for this we have the standard
orthonormal basis consisting of the following elements
$$\frac 1{\sqrt 2}\begin{pmatrix}Y_{rs} & 0 \\
0 & Y_{rs}\end{pmatrix},\frac 1{\sqrt 2}\begin{pmatrix}X_{rs} & 0 \\
0 & -X_{rs}\end{pmatrix},\frac 1{\sqrt 2}\begin{pmatrix}D_t &0 \\
0 & -D_t\end{pmatrix},$$

$$\frac 1{\sqrt 2}\begin{pmatrix}0 & X_{rs} \\
X_{rs} & 0\end{pmatrix},
\frac 1{\sqrt 2}\begin{pmatrix}0 & X_{rs} \\
-X_{rs} & 0\end{pmatrix},\frac 1{\sqrt 2}\begin{pmatrix}0 &D_t \\
D_t & 0\end{pmatrix},
\frac 1{\sqrt 2}\begin{pmatrix}0 &D_t \\
-D_t& 0\end{pmatrix}$$ where $1\le r<s\le n$ and $1\le t\le n$.

\begin{lemma}\label{lemm:Sp(n,R)}
For $1\le i,j\le n$ let $x_{ij},y_{ij},z_{ij},w_{ij}:\SpR n\to\cn$
be the real valued coordinate functions given by
$$x_{ij}:g\mapsto e_i\cdot g\cdot e_j^t,\ \
y_{ij}:g\mapsto e_i\cdot g\cdot e_{n+j}^t, $$
$$z_{ij}:g\mapsto e_{n+i}\cdot g\cdot e_j^t,\ \
w_{ij}:g\mapsto e_{n+i}\cdot g\cdot e_{n+j}^t$$ where $\{e_1,\dots
,e_{2n}\}$ is the canonical orthonormal basis for $\rn^{2n}$. Then
the following relations hold
$$\tau(x_{ij})=\frac 12 x_{ij}, \ \
\tau(y_{ij})=\frac 12 y_{ij}, \ \ \tau(z_{ij})=\frac 12 z_{ij}, \
\ \tau(w_{ij})=\frac 12 w_{ij},$$
$$\kappa(x_{ij}, x_{kl})= \frac{1}{2} \bigl[y_{il} y_{kj} +
\delta _{jl}\cdot\sum _{t=1}^n( x_{it} x_{kt}+y_{it}
y_{kt})\bigr],$$
$$\kappa(x_{ij},y_{kl})=-\frac{1}{2}x_{il}y_{kj},$$
$$\kappa(x_{ij}, z_{kl})= \frac{1}{2} \bigl[y_{il} w_{kj} +
\delta _{jl}\cdot\sum _{t=1}^n( x_{it} z_{kt}+y_{it}
w_{kt})\bigr],$$
$$\kappa(x_{ij},w_{kl})=-\frac{1}{2}x_{il}w_{kj},$$
$$\kappa(y_{ij}, y_{kl})= \frac{1}{2} \bigl[x_{il} x_{kj} +
\delta _{jl}\cdot\sum _{t=1}^n( x_{it} x_{kt}+y_{it}
y_{kt})\bigr],$$
$$ \kappa(y_{ij},z_{kl})=-\frac{1}{2}y_{il}z_{kj},$$
$$\kappa(y_{ij}, w_{kl})= \frac{1}{2} \bigl[x_{il} z_{kj} +
\delta _{jl}\cdot\sum _{t=1}^n( x_{it} z_{kt}+y_{it}
w_{kt})\bigr],$$
$$\kappa(z_{ij},z_{kl})= \frac{1}{2} \bigl[w_{il} w_{kj} +
\delta _{jl}\cdot\sum _{t=1}^n( z_{it} z_{kt}+w_{it}
w_{kt})\bigr],$$
$$\kappa(z_{ij},w_{kl})=-\frac{1}{2}z_{il}w_{kj},$$
$$\kappa(w_{ij}, w_{kl})= \frac{1}{2} \bigl[z_{il} z_{kj} +
\delta _{jl}\cdot\sum _{t=1}^n( z_{it} z_{kt}+w_{it}
w_{kt})\bigr].$$
\end{lemma}

\begin{proof}  The method of proof is exactly the same as for
Lemma \ref{lemm:special}.
\end{proof}

\begin{theorem}\label{theo:Sp(n,R)-1}
Let $v\in\cn^n$ be a non-zero element, then the $2n$-dimensional
vector space
$$\E (v)=\{\phi_{ab}(g)=\trace (a^t v (x+iy)^t+b^t v(z+ i w)^t)
\ |\ a,b\in \cn ^n \},$$ is an eigenfamily on $\SpR n$.
\end{theorem}

\begin{theorem}\label{theo:Sp(n,R)-2}
For non-zero elements $a,b\in \cn ^n$ the vector space
$$\E (a,b)=\{\phi_v(g)=
\trace (a^t v (x+iy)^t+b^t v(z+ i w)^t)\ |\ v\in \cn ^n \},$$ is
an eigenfamily on $\SpR n$.
\end{theorem}

\begin{proof} Employing the results of Lemma \ref{lemm:Sp(n,R)}
it is easily seen that if $P$ and $Q$ both belong to one of the sets
described above then $2\kappa(P,Q)+PQ=0$.
\end{proof}

\section{Bi-Eigenfamilies}

The notion of an eigenfamily has in earlier sections turned out to
be useful for the construction of harmonic morphisms from the
Riemannian non-compact semisimple Lie groups $\SLR n$, $\SUs {2n}$
and $\SpR n$.  In order to solve the problem in the cases of
$$\SOs{2n}´,\ \ \SOO pq,\ \ \SUU pq\ \ \text{and}\ \ \Spp pq$$ we
now generalize to the notion of a bi-eigenfamily.

\begin{definition}\label{defi:bi-eigen}
Let $(M,g)$ be a semi-Riemannian manifold and $$\E_1
=\{\phi_i:M\to\cn\ |\ i\in I_1\},\ \ \E_2 =\{\psi_i:M\to\cn\ |\
i\in I_2\}$$ be two non-empty eigenfamilies on $M$ i.e. there
exist complex numbers $\lambda_{1}$, $\mu_1$, $\lambda_{2}$,
$\mu_2$ such that
$$\tau(\phi_1)=\lambda_{1}\phi_1,\ \
\kappa (\phi_1,\phi_2)=\mu_1\phi_1\phi_2,$$
$$\tau(\psi_1)=\lambda_{2}\psi_1,\ \
\kappa (\psi_1,\psi_2)=\mu_2\psi_1\psi_2$$ for all
$\phi_1,\phi_2\in\E_1$ and $\psi_1,\psi_2\in\E_2$. The union $$\B
=\E_1 \cup \E_2$$ is said to be a {\it bi-eigenfamily} on $M$ if
there exists a complex number $\mu$ such that
$$\kappa (\phi,\psi)=\mu\phi\psi$$ for all $\phi\in\E_1$
and $\psi\in\E_2$.
\end{definition}

\begin{definition}\label{defi:bi-homo}
Let $m,n,d_1,d_2$ be non-negative integers.  Then a polynomial map
$P:\cn^m\times\cn^n\to\cn$ is said to be {\it bi-homogeneous} of
{\it bi-degree} $d=(d_1,d_2)$ if for all $\lambda,\mu\in\cn$
$$P(\lambda z,\mu w)=\lambda^{d_1}\mu^{d_2} P(z,w).$$
\end{definition}

The following result generalizes that of Theorem
\ref{theo:rational} and shows that bi-eigenfamilies are useful
ingredients for the construction of harmonic morphisms.

\begin{theorem}\label{theo:bi-rational}
Let $(M,g)$ be a semi-Riemannian manifold and $\B =\E_1 \cup \E_2$
be a finite bi-eigenfamily on $M$ where
$$\E_1=\{\phi_1,\dots,\phi_m\}\ \ \text{and}\ \ \E_2
=\{\psi_1,\dots,\psi_n\}.$$ If $P,Q:\cn^m\times\cn^n\to\cn$ are
linearly independent bi-homogeneous polynomials of the same
bi-degree then the quotient
$$\frac{P(\phi_1,\ldots,\phi_m,\psi_1,\ldots,\psi_n)}
{Q(\phi_1,\ldots,\phi_m,\psi_1, \ldots,\psi_n)}$$ is a
non-constant harmonic morphism on the open and dense subset
$$\{p\in M\ | \ Q(\phi_1(p), \ldots, \phi_m(p), \psi_1(p),
\ldots, \psi_n(p))\neq 0\}.$$
\end{theorem}

\begin{proof} A proof of Theorem \ref{theo:bi-rational} can be
found in Appendix \ref{app:general}.\end{proof}

\section{The Riemannian Lie group $\SOs{2n}$}

In this section we construct eigenfamilies and bi-eigenfamilies of
complex valued functions on the non-compact semisimple Lie group
$$\SOs{2n}=\{q\in\SUU nn\ |\ q\cdot I_{nn}\cdot J_{n}\cdot
q^t=I_{nn}\cdot J_{n}\}$$ where
$$\UU nn =\{z\in\GLC{2n}\ |\ z\cdot I_{n,n}\cdot z^*=I_{n,n}\}$$
and $$I_{nn}=\begin{pmatrix} -I_n& 0\\0 & I_n\end{pmatrix},\ \
J_n=\begin{pmatrix}0 & I_n \\-I_n & 0\end{pmatrix}.$$

The Lie algebra $\sos{2n}$ of $\SOs{2n}$ is given
by $$\sos{2n}=\{\begin{pmatrix}Z & W \\
-\bar{W} & \bar{Z}\end{pmatrix}\in\cn^{2n\times 2n}|\ Z+Z^*=0\
\text{and}\ W+W^t=0\}$$  and for this we have the standard
orthonormal basis consisting of the following matrices
$$\frac 1{\sqrt 2}\begin{pmatrix}Y_{rs} & 0 \\
0 & Y_{rs}\end{pmatrix},\frac 1{\sqrt 2}\begin{pmatrix}iX_{rs} & 0 \\
0 & -iX_{rs}\end{pmatrix},\frac 1{\sqrt 2}\begin{pmatrix}iD_t & 0 \\
0 & -iD_t\end{pmatrix},$$

$$\frac 1{\sqrt 2}\begin{pmatrix}0 & Y_{rs} \\
-Y_{rs} & 0\end{pmatrix},\frac 1{\sqrt 2}\begin{pmatrix}0 & iY_{rs} \\
iY_{rs} & 0\end{pmatrix}$$ where $1\le r<s\le n$ and $1\le t\le
n$.

\begin{lemma}\label{lemm:SO*(2n)}
For $1\le i,j\le n$ let $z_{ij},w_{ij}:\SOs{2n}\to\cn$ be the
complex valued coordinate functions given by
$$z_{ij}:q\mapsto e_i\cdot q\cdot e_j^t,\ \
w_{ij}:q\mapsto e_i\cdot q\cdot e_{n+j}^t$$ where $\{e_1,\dots
,e_{2n}\}$ is the canonical orthonormal basis for $\cn^{2n}$. Then
the following relations hold
$$\tau(z_{ij})=-\frac{1}{2}z_{ij},\ \ \tau(w_{ij})=-\frac{1}{2}w_{ij},$$
$$\kappa(z_{ij}, z_{kl})= -\frac{1}{2}z_{il}z_{kj}, \ \ \kappa(w_{ij},
w_{kl})= -\frac{1}{2}w_{il}w_{kj},$$
$$\kappa(z_{ij}, w_{kl})=\frac{1}{2} \big[z_{kj}w_{il}+
\delta _{jl}\cdot\sum _{t=1}^n(z_{it} w_{kt}-w_{it}
z_{kt})\big].$$
\end{lemma}

\begin{proof}  The method of proof is exactly the same as for
Lemma \ref{lemm:special}.
\end{proof}

\begin{theorem}\label{theo:SO*(2n)}
Let $v$ be a non-zero element of $\cn ^n$. Then the vector spaces
$$\E_1(v)=\{\phi_{a}:\SOs{2n}\to\cn\ |\ \phi_{a}(z+jw)=
\trace (v^ta z^t ),\ a\in \cn ^n\} $$ and
$$\E_2(v)=\{\psi_{c}:\SOs{2n}\to\cn\ |\ \psi_{c}(z+jw)=
\trace (v^tc w^t),\ c\in \cn^n\} $$ are eigenfamilies on
$\SOs{2n}$ and their union $$\B =\E_1(v)\cup\E_2(v)$$ is a
bi-eigenfamily.
\end{theorem}

\begin{proof}  If $\phi_a,\phi_b\in\E_1(v)$ and
$\psi_c,\psi_d\in\E_2(v)$ then
it is an immediate consequence of Lemma \ref{lemm:SO*(2n)} that
$$\tau(\phi_a)=-\frac 12\phi_a,\ \ \tau(\psi_c)=-\frac 12\psi_c,$$
$$\kappa (\phi_a,\phi_b)=-\frac 12\phi_a\phi_b,\ \
\kappa (\phi_b,\psi_c)=\frac 12\phi_b\psi_c,\ \ \kappa
(\psi_c,\psi_d)=-\frac 12\psi_c\psi_d.$$  This proves the
statement of Theorem \ref{theo:SO*(2n)}.
\end{proof}

\begin{theorem}\label{theo:SO*(2n)-2}
Let $a$ be a non-zero element of $\cn ^n$. Then the vector spaces
$$\E_1(a)=\{\phi_{u}:\SOs{2n}\to\cn\ |\ \phi_{u}(z+jw)=
\trace (u^ta z^t ),\ u\in \cn ^n\} $$ and
$$\E_2(a)=\{\psi_{v}:\SOs{2n}\to\cn\ |\ \psi_{v}(z+jw)=
\trace (v^ta w^t),\ v\in \cn^n\} $$ are eigenfamilies on
$\SOs{2n}$.
\end{theorem}

\section{Notation}

In this section we introduce the necessary notation for dealing
with the important non-compact cases $\SOO pq$, $\SUU pq$ and
$\Spp pq$.  Further we state a few matrix identities useful for
those readers interested in the details of our calculations.

For a positve integer $n$ let $p,q$ be non-negative integers
such that $n=p+q$.  Further let the index sets
$\Delta_1,\Delta_2,\Lambda_1,\Lambda_2$ be defined by
$$\Delta _1=\{1, \ldots, p\},\ \ \Delta _2=\{p+1, \ldots, n\},$$
$$\Lambda _1=\{(r,s) \ | \ 1\le r<s\le p\ \ \text{or}
\ \  p+1\le r<s \le n \},$$
$$\Lambda _2=\{(r,s)\ | \ 1\le r\le p\ \ \text{and}\ \ p+1\le s \le n \}$$
and the matrices $I_{pq},J_n$ by
$$I_{pq}=\begin{pmatrix} -I_p& 0\\0 & I_q\end{pmatrix}\ \ \text{and}\ \
J_n=\begin{pmatrix}0 & I_n \\-I_n & 0\end{pmatrix}.$$ With the
above notation we have the following easily verified matrix
identities
$$\sum_{(r,s) \in \Lambda _1}X_{rs}^2 -\sum_{{(r,s) \in
\Lambda _2}}X_{rs}^2=-\frac{1}{2} (I_n +(p-q)I_{pq}),$$

$$\sum_{(r,s) \in \Lambda _1}Y_{rs}^2 -\sum_{{(r,s) \in
\Lambda _2}}Y_{rs}^2=\frac{1}{2} (I_n +(p-q)I_{pq}),$$

\begin{eqnarray*}
&&\sum_{(r,s) \in \Lambda _1}X_{rs}E_{jl}X_{rs}^t
-\sum_{(r,s) \in \Lambda _2}X_{rs}E_{jl}X_{rs}^t\\
&&\qquad\qquad\qquad=\frac{1}{2}((-1)^{\chi (j)+\chi (l)+\delta
_{jl}}E_{lj} +\delta_{jl} (-1)^{\chi   (j)}I_{pq}),
\end{eqnarray*}

\begin{eqnarray*}
&&\sum_{(r,s) \in \Lambda _1}Y_{rs}E_{jl}Y_{rs}^t
-\sum_{(r,s) \in \Lambda _2} Y_{rs}E_{jl}Y_{rs}^t\\
&&\qquad\qquad\qquad=\frac{1}{2}((-1)^{\chi (j)+\chi (l)+1}E_{lj}
+\delta_{jl} (-1)^{\chi (j)}I_{pq}).
\end{eqnarray*}

\begin{eqnarray*}
&&\sum_{(r,s) \in \Lambda _1}X_{rs}E_{jl}X_{rs}^t
 +\sum_{(r,s) \in \Lambda _2}Y_{rs}E_{jl}Y_{rs}^t
 +\sum_{t=1}^nD_{t}E_{jl}D_{t}^t\\
&&\qquad\qquad\qquad=\frac{1}{2}((-1)^{\chi (j)+\chi (l)+1}E_{lj}
+\delta_{jl} (-1)^{\chi (j)}I_{n}).
\end{eqnarray*}
Here $\chi:\zn\to\{0,1\}$ is the characteristic function for the
set $\Delta_1$ given by
$$\chi (i)=\begin{cases}
1,& i \in \Delta _1,\\
0,& i \notin \Delta _1.
\end{cases}$$
In the sequel we shall use the notation $$\cn^p_1
=\{(z,w)\in\cn^p\times\cn^q\;|\;w=0\},\ \
\cn^q_2=\{(z,w)\in\cn^p\times\cn^q\;|\;z=0\}.$$

\section{The Riemannian Lie group $\SUU pq$}

In this section we construct eigenfamilies on the non-compact Lie
group $\UU pq$. These can be used to construct a variety of local
harmonic morphisms on the semisimple
$$\SUU pq=\{z\in\SLC{p+q}\ |\ z\cdot I_{p,q}\cdot z^*=I_{p,q}\}.$$
The Lie group $\UU pq$ is the non-compact subgroup of $\GLC n$
given by
$$\UU pq =\{z\in\GLC{p+q}\ |\ z\cdot I_{p,q}\cdot z^*=I_{p,q}\}.$$
The circle group $S^1=\{w\in\cn\ | \ |w|=1\}$ acts on $\UU pq$ by
multiplication $(w,z)\mapsto wz$ and the orbit space of this
action is the group $\SUU pq$. The Lie algebra $\uu pq$ of $\UU
pq$ is given by
$$\uu pq=\{\begin{pmatrix}Z & M \\
M^* & W\end{pmatrix}\in\glc {p+q} |\ Z+Z^t=0\ \text{and}\
W+W^t=0\}$$ with the canonical orthonormal basis
$$\{iX_{rs}, Y_{rs}\ |\ (r,s)\in\Lambda_1\}\cup
   \{X_{rs},iY_{rs}\ |\ (r,s)\in\Lambda_2\}\cup
   \{iD_t\ |\ t=1,\dots ,n\}.$$

\begin{lemma}\label{lemm:complex}
For $1\le k,l\le n$ let $z_{kl}:\UU pq\to\cn$ be the complex
valued coordinate functions given by
$$z_{kl}:z\mapsto e_k\cdot z\cdot e_l^t$$ where $\{e_1,\dots
,e_n\}$ is the canonical orthonormal basis for $\cn^n$. Then the
following relations hold
\begin{eqnarray*}
\tau (z_{ij})&=&(-1)^{\chi (j)} (p-q) z_{ij},\\
\kappa(z_{ij},z_{kl})&=&-(-1)^{\chi (j)+\chi (l)}z_{il} z_{kj}.
\end{eqnarray*}
\end{lemma}

\begin{proof}
It follows directly from the definition of the functions $z_{ij}$
that if $Z$ is an element of the Lie algebra $\uu pq$ then the
first and second order derivatives satisfy
$$Z(z_{ij}):z\mapsto e_i\cdot z\cdot Z\cdot e_j^t\ \ \text{and}
\ \  Z^2(z_{ij}):z\mapsto e_i\cdot z\cdot Z^2\cdot e_j^t.$$
Employing the above mentioned matrix identities we then yield
\begin{eqnarray*}
\tau(z_{ij}) &=&\sum_{(r,s)\in\Lambda _1}Y_{rs}^2(z_{ij})
-\sum_{(r,s)\in\Lambda_2}Y_{rs}^2
(z_{ij})-\sum_{(r,s)\in\Lambda_1}X_{rs}^2(z_{ij})\\
&&\quad +\sum_{(r,s)\in\Lambda_2}X_{rs}^2 (z_{ij})-\sum_{t=1} ^n D_t^2(z_{ij})\\
&=&\sum_{(r,s)\in\Lambda_1}e_izY_{rs}^2e_j ^t -\sum_{(r,s) \in
\Lambda _2}e_izY_{rs}^2e_j ^t
-\sum_{(r,s) \in \Lambda _1} e_i z X_{rs}^2e_j ^t \\
& &\quad +\sum_{(r,s) \in \Lambda _2}e_i z X_{rs}^2e_j
^t -\sum_{t=1} ^n e_i z D_t^2e_j ^t\\
&=&(p-q)
e_i z I_{pq}e_j ^t\\
&=&(-1)^{\chi(j)}(p-q)z_{ij}.
\end{eqnarray*}

\begin{eqnarray*}
\kappa(z_{ij},z_{kl}) &=&\sum_{(r,s) \in \Lambda _1}e_iz
Y_{rs}e_j^t e_l Y_{rs}^t z^te_k^t-\sum_{(r,s) \in \Lambda
_2}e_izY_{rs}e_j^te_l
Y_{rs}^tz^te_k^t\\
& &-\sum_{(r,s) \in \Lambda _1} e_izX_{rs}e_j^t e_l X_{rs}^t
z^te_k^t+\sum_{(r,s) \in \Lambda _2}e_izX_{rs}e_j^te_l
X_{rs}^tz^te_k^t\\
& &-\sum_{t=1}^n e_izD_te_j^te_lD_t^tz^te_k^t\\
&=&e_i z \Bigl(\sum_{(r,s) \in \Lambda _1}Y_{rs}E_{jl}
Y_{rs}^t-\sum_{(r,s) \in \Lambda _2}Y_{rs}E_{jl}Y_{rs}^t\Bigr)
z^t e_k ^t\\
&&-e_i z \Bigl(\sum_{(r,s) \in \Lambda
_1}X_{rs}E_{jl}X_{rs}^t-\sum_{(r,s) \in \Lambda _2}
X_{rs}E_{jl}X_{rs}^t\Bigr)z^te_k ^t\\
&&\quad -e_i z
\sum_{t=1}^n D_t E_{jl}D_t^tz^te_k ^t\\
&=&-(-1)^{\chi (j)+\chi (l)}z_{il} z_{kj}.
\end{eqnarray*}
\end{proof}

Let $P,Q:\UU pq\to\cn$ be homogeneous polynomials of the
coordinate functions $z_{ij}:\UU pq\to\cn$ of degree one i.e of
the form
$$P=\trace (A\cdot z^t)=\sum_{i,j=1}^n a_{ij}z_{ij}\ \ \text{and}
\ \ Q=\trace (B\cdot z^t)=\sum_{k,l=1}^n b_{kl}z_{kl}$$ for some
$A,B\in\cn^{n\times n}$.  As a direct consequence of Lemma
\ref{lemm:complex} we then yield
$$PQ+\kappa(P,Q)=\sum_{i,j,k,l=1}^n(a_{ij}b_{kl}-
(-1)^{\chi (j)+\chi (l)}a_{kj}b_{il})z_{ij}z_{kl}.$$

\begin{theorem}\label{theo:complex-1}
Let $v$ be a non-zero element of $\cn ^n$. Then the vector spaces
$$\E_1(v)=\{\phi_a:\UU pq\to\cn\ |\ \phi_a(z)=\trace (v^taz^t),\
a\in \cn^p_1\}$$ and
$$\E_2 (v)=\{\phi_a:\UU pq\to\cn\ |\ \phi_a(z)=\trace (v^taz^t),\
a\in \cn^q_2\}$$ are eigenfamilies on $\UU pq$ and their union
$$\B=\E_1(v)\cup\E_2(v)$$ is a bi-eigenfamily.
\end{theorem}

\begin{proof} Assume that $a,b\in\cn^p_1$ and define $A=v^ta$ and
$B=v^tb$. By construction any two columns of the matrices $A$ and
$B$ are linearly dependent.  This means that for all $1\le
i,j,k,l\le n$
$$\det\begin{pmatrix} a_{ij}& b_{il}
\\a_{kj} & b_{kl}\end{pmatrix} =(a_{ij}b_{kl}-a_{kj}b_{il})=0$$
so $P^2+\kappa (P,P)=0$, $PQ+\kappa (P,Q)=0$ and $Q^2+\kappa
(Q,Q)=0$. The fact that $\E_1(v)$ is an eigenfamily is now a
direct consequence of Lemma \ref{lemm:complex}.  A similar
argument shows that $\E_2(v)$ is also an eigenfamily.  Note that
$\kappa(P,Q)=PQ$ for any $P \in \E_1 (v)$, $Q \in \E_2 (v)$.  This
shows that $\B$ is a bi-eigenfamily on $\UU pq$.
\end{proof}

\begin{theorem}\label{theo:complex-2}
Let $u,v$ be two non-zero elements of $\cn^p_1$ and $\cn^q_2$,
respectively. Then the sets
$$\E_1(u)=\{\phi_a:\UU pq\to\cn\ |\ \phi_a(z)=\trace (a^tuz^t),\
a\in \cn^n\},$$
$$\E_2(v)=\{\phi_a:\UU pq\to\cn\ |\ \phi_a(z)=\trace (a^tvz^t),\
a\in \cn^n\}$$ are eigenfamilies on $\UU pq$.
\end{theorem}

\begin{proof}  The argument is similar to that used to prove the
first part of Theorem \ref{theo:complex-1}.
\end{proof}

It should be noted that the local harmonic morphisms on $\UU pq$
that we obtain by applying Theorem \ref{theo:bi-rational} are
invariant under the cirle action and hence induce local harmonic
morphisms on the semisimple Lie group $\SUU pq$.

\section{The Riemannian Lie group $\SOO pq$}

In this section we construct eigenfamilies on the non-compact
semisimple Lie group $$\SOO pq=\{x\in\SLR{p+q}\ |\ x\cdot
I_{p,q}\cdot x^t=I_{p,q}\}.$$

The Lie algebra $\soo pq$ of $\SOO pq$ is given by
$$\soo pq=\{\begin{pmatrix}A & M \\
M^t & B\end{pmatrix}\in\slr {p+q} |\ A+A^t=0\ \text{and}\
B+B^t=0\}.$$  For this we have the canonical orthonormal basis
$$\{Y_{rs}\ |\ (r,s)\in\Lambda_1\}\cup\{X_{rs}\ |\ (r,s)\in\Lambda_2\}.$$

\begin{lemma}\label{lemm:real}
For $1\le i,j\le n$ let $x_{ij}:\SOO pq\to\rn$ be the real valued
coordinate functions given by
$$x_{ij}:x\mapsto e_i\cdot x\cdot e_j^t$$ where $\{e_1,\dots
,e_n\}$ is the canonical orthonormal basis for $\rn^n$. Then the
following relations hold

\begin{eqnarray*}
\tau (x_{ij})&=&\frac 12(1+(-1)^{\chi(j)}(p-q))x_{ij},\\
\kappa(x_{ij},x_{kl})&=&\frac{1}{2}\bigl((-1)^{\chi(j)+\chi(l)+1}x_{il}x_{kj}\\
& &\quad -\delta_{jl} (-1)^{\chi (j)} [ \sum _{t\in\Delta _1 }
x_{it} x_{kt} -\sum _{t\in\Delta _2 } x_{it}x_{kt} ]\bigr).
\end{eqnarray*}
\end{lemma}

\begin{proof}
The proof is similar to that of Lemma \ref{lemm:complex}.
\end{proof}

Let $P,Q:\SOO pq\to\cn$ be homogeneous polynomials of the
coordinate functions $x_{ij}:\SOO pq\to\cn$ of degree one i.e. of
the form
$$P=\trace (A\cdot x^t)=\sum_{i,j=1}^n a_{ij}x_{ij}\ \ \text{and}
\ \ Q=\trace (B\cdot x^t)=\sum_{k,l=1}^n b_{kl}x_{kl}$$ for some
$A,B\in\cn^{n\times n}$.  As a direct consequence of Lemma
\ref{lemm:real} we then yield
\begin{eqnarray*}
PQ+2\kappa(P,Q)&=&\sum_{i,j,k,l=1}^na_{ij}b_{kl}x_{ij}x_{kl}
+2\sum_{i,j,k,l=1}^na_{ij}b_{kl}\kappa(x_{ij},x_{kl})\\
&=&\sum_{i,j,k,l=1}^na_{ij}b_{kl}x_{ij}x_{kl}
-(-1)^{\chi(j)+\chi(l)}\sum_{i,j,k,l=1}^na_{ij}b_{kl}x_{il}x_{kj}\\
&&\quad -\delta_{jl} (-1)^{\chi (j)}\sum_{i,j,k,=1}^na_{ij}b_{kj}
[\sum_{t=1}^px_{it}x_{kt}-\sum_{t=p+1}^{p+q}x_{it}x_{kt}]\\
&=&\sum_{i,j,k,l=1}^n(a_{ij}b_{kl}-(-1)^{\chi(j)+\chi(l)}a_{kj}b_{il})x_{ij}x_{kl}\\
&&\quad -\delta_{jl} (-1)^{\chi (j)}\sum_{i,j,k,=1}^na_{ij}b_{kj}
[\sum_{t=1}^px_{it}x_{kt}-\sum_{t=p+1}^{p+q}x_{it}x_{kt}]\\
\end{eqnarray*}

Comparing coefficients in the case when the integer
$\chi(j)+\chi(l)$ is even we see that $PQ+2\kappa (P,Q)=0$ if for
all $1\le i,k\le n$ $$\sum_{j=1}^na_{ij}b_{kj}=0$$ and
$$\det\begin{pmatrix} a_{ij}& b_{il} \\a_{kj} &
b_{kl}\end{pmatrix} =(a_{ij}b_{kl}-a_{kj}b_{il})=0.$$ for all
$1\le i,j,k,l\le n$.

\begin{theorem}\label{theo:real-1}
Let $u$ be a non-zero element of $\cn ^n$ and $V_1$, $V_2$ be
maximal isotropic subspaces of $\cn^p_1$ and $\cn^q_2$,
respectively.  Then the vector spaces
$$\E_{V_1}(u)=\{\phi_a:\SOO pq\to\cn\ |\ \phi_a(x)=\trace (u^tax^t),\
a\in V_1\}$$ and
$$\E_{V_2}(u)=\{\phi_a:\SOO pq\to\cn\ |\ \phi_a(x)=\trace (u^tax^t),\
a\in V_2\}$$ are eigenfamilies on $\SOO pq$ and their union
$$\B=\E_{V_1}(u)\cup\E_{V_2}(u)$$ is a bi-eigenfamily.
\end{theorem}

\begin{proof}  Assume that $a,b\in V_1$ and define $A=u^ta$ and
$B=u^tb$.  By construction any two columns of the matrices $A$ and
$B$ are linearly dependent.  This means that for all $1\le
i,j,k,l\le n$
$$\det\begin{pmatrix} a_{ij}& b_{il}
\\a_{kj} & b_{kl}\end{pmatrix} =(a_{ij}b_{kl}-a_{kj}b_{il})=0.$$
Furthermore we have
$$PQ+2\kappa (P,Q)=\sum_{j=1}^na_{ij}b_{kj}=u_iu_k(a,b)=0.$$
Hence $P^2+2\kappa (P,P)=0$, $PQ+2\kappa (P,Q)=0$, $Q^2+2\kappa
(Q,Q)=0$. It follows directly from Lemma \ref{lemm:real} that
$\E_{V_1}(u)$ is an eigenfamily.  A similar argument shows that
even $\E_{V_2}(u)$ is an eigenfamily on $\SOO pq$. It is easy to
see that $2\kappa(P,Q)-PQ=0$ for any $P \in \E_{V_1}(u)$, $Q \in
\E_{V_2}(u)$. This shows that $\B$ is a bi-eigenfamily on $\SOO
pq$.
\end{proof}

\begin{theorem}\label{theo:real-2}
Let $u\in\cn^p_1$ and $v\in\cn^q_2$ be two non-zero isotropic
elements of $\cn^n$ i.e. $(u,u)=(v,v)=0$. Then the sets
$$\E_1 (u)=\{\phi_a:\SOO pq\to\cn\ |\ \phi_a(x)=\trace (a^tux^t),\
a\in \cn^n\}$$ and
$$\E_2(v)=\{\phi_a:\SOO pq\to\cn\ |\ \phi_a(x)=\trace (a^tvx^t),\
a\in \cn^n\}$$ are eigenfamilies on $\SOO pq$.
\end{theorem}

\begin{proof}  The argument is similar to that used to prove the
first part of Theorem \ref{theo:real-1}.
\end{proof}

\section{The Riemannian Lie group $\Spp pq$}

In this section we construct eigenfamilies on the non-compact
semisimple Lie group
$$\Spp pq =\{g\in\GLH{p+q}\ |\ g\cdot I_{p,q}\cdot g^*=I_{p,q}\}.$$
Using the standard representation of the quaternions $\hn$ in
$\cn^{2\times 2}$
$$(z+jw)\mapsto g=\begin{pmatrix}z & w \\ -\bar w & \bar
z\end{pmatrix}$$ it is easily seen that the Lie algebra $\spp pq$
of $\Spp pq$ is the set of matrices
$$\begin{pmatrix}
Z_{11}   & Z_{12} & W_{11}   & W_{12}\\
\bar Z^t_{12} & Z_{22} & -W_{12}^t & W_{22}\\
-\bar W_{11} & -\bar W_{12} & \bar Z_{11} & \bar Z_{12}\\
\bar W^t_{12} & -\bar W_{22} & Z^t_{12} & \bar Z_{22}\end{pmatrix}
\in\cn^{2n\times 2n}$$ where $Z_{11}+\bar Z^t_{11}=0$,
$Z_{22}+\bar Z^t_{22}=0$, $W_{11}=W_{11}^t$ and $W_{22}=W^t_{22}$.
For the Lie algebra we have the standard orthonormal basis being
the union of the following sets
$$\{\frac 1{\sqrt 2}\begin{pmatrix}Y_{rs} & 0 \\
0 & Y_{rs}\end{pmatrix},\frac 1{\sqrt 2}\begin{pmatrix}iX_{rs} & 0 \\
0 & -iX_{rs}\end{pmatrix}\ |\ (r,s)\in\Lambda_1\},$$

$$\{\frac 1{\sqrt 2}\begin{pmatrix}X_{rs} & 0 \\
0 & X_{rs}\end{pmatrix},\frac 1{\sqrt 2}\begin{pmatrix}iY_{rs} & 0 \\
0 & -iY_{rs}\end{pmatrix}\ |\ (r,s)\in\Lambda_2\},$$

$$\{\frac 1{\sqrt 2}\begin{pmatrix}0 & X_{rs} \\
-X_{rs} & 0\end{pmatrix},\frac 1{\sqrt 2}\begin{pmatrix}0 & iX_{rs} \\
iX_{rs} & 0\end{pmatrix}\ |\ (r,s)\in\Lambda_1\},$$

$$\{\frac 1{\sqrt 2}\begin{pmatrix}0 & Y_{rs} \\
-Y_{rs} & 0\end{pmatrix},\frac 1{\sqrt 2}\begin{pmatrix}0 & iY_{rs} \\
iY_{rs} & 0\end{pmatrix}\ |\ (r,s)\in\Lambda_2\},$$

$$\{\frac 1{\sqrt 2}\begin{pmatrix}iD_{t} & 0 \\
0 & -iD_{t}\end{pmatrix},\frac 1{\sqrt 2}\begin{pmatrix}0 & D_{t}  \\
-D_{t} & 0\end{pmatrix},\frac 1{\sqrt 2}\begin{pmatrix}0 & iD_{t}  \\
iD_{t} & 0\end{pmatrix}\ |\ 1\le t\le n\}.$$

\begin{lemma}\label{lemm:quaternionic}
For $1\le i,j\le n=p+q$ let $z_{ij},w_{ij}:\Spp pq\to\cn$ be the
complex valued coordinate functions given by
$$z_{ij}:g\mapsto e_i\cdot g\cdot e_j^t,\ \
w_{ij}:g\mapsto e_i\cdot g\cdot e_{n+j}^t$$ where $\{e_1,\dots
,e_{2n}\}$ is the canonical orthonormal basis for $\cn^{2n}$. Then
the following relations hold
$$\tau(z_{ij})=-\frac 12\bigl[(-1)^{\chi(l)}2(q-p)+1\bigr]z_{ij}.$$
$$\tau(w_{ij})=-\frac 12\bigl[(-1)^{\chi(l)}2(q-p)+1\bigr]w_{ij}.$$
$$\kappa(z_{ij}, z_{kl})= -\frac{1}{2}(-1)^{\chi
(j)+\chi(l)}z_{il}z_{kj}$$
$$\kappa(w_{ij}, w_{kl})= -\frac{1}{2}(-1)^{\chi
(j)+\chi(l)}w_{il}w_{kj}$$
$$\kappa(z_{ij}, w_{kl})= -\frac{1}{2}\bigl[(-1)^{\chi
(j)+\chi(l)}w_{il}z_{kj}-\delta _{jl}\cdot\sum _{t=1}^n( z_{it}
w_{kt}-w_{it} z_{kt})\bigr]$$
\end{lemma}

\begin{proof}  The proof is similar to that of Lemma \ref{lemm:complex}
but more involved.
\end{proof}

\begin{theorem}\label{theo:quaternionic-1}
Let $v$ be a non-zero element of $\cn ^n$. Then the vector spaces
$$\E_1 (v)=\{\phi_{ab}:\Spp pq\to\cn\ |\ \phi_{ab}(z+jw)=
\trace (v^ta z^t + v^tbw^t),\ a,b\in \cn^p_1\} $$ and
$$\E_2 (v)=\{\phi_{ab}:\Spp pq\to\cn\ |\ \phi_{ab}(z+jw)=
\trace (v^t a z^t + v^tbw^t),\ a,b\in \cn^p_2\} $$ are eigenfamilies
on $\Sp n$ and their union
$$\B=\E_1(v)\cup\E_2(v)$$ is a bi-eigenfamily.
\end{theorem}

\begin{proof} Let $a,b,c,d$ be arbitrary elements of $\cn^p_1$ and
define the complex valued functions $P,Q:\Spp pq\to\cn$ by
$$P=\trace (v^taz^t+v^tbw^t)\ \
\text{and}\ \ Q=\trace (v^tcz^t+v^tdw^t).$$ Then a simple
calculation shows that
$$PQ+2\kappa(P,Q)=[(a,d)-(b,c)]
\sum_{i,k,t=1}^n v_i v_k (z_{it}w_{kt}-w_{it}z_{kt})=0.$$
Automatically we also get $P^2+2\kappa(P,P)=0$ and
$Q^2+2\kappa(Q,Q)=0$. It now follows directly from Lemma
 \ref{lemm:real} that $\E_1 (v)$ is an eigenfamily.  A similar
argument shows that $\E_2 (v)$ is an eigenfamily on $\Spp pq$. It
is easy to see that $2\kappa(P,Q)-PQ=0$ for any $P \in \E_1 (v)$,
$Q \in \E_2 (v)$. This shows that $\B$ is a bi-eigenfamily on
$\Spp pq$.
\end{proof}

\section{The Duality}

In this section we show how a real analytic bi-eigenfamily $\B$ on
a Riemannian non-compact semisimple Lie group $G$ gives rise to a
real-analytic bi-eigenfamily $\B^*$ on its semi-Riemannian compact
dual $U$ and vice versa.  The method of proof is borrowed from a
related duality principle for harmonic morphisms from Riemannian
symmetric spaces, see \cite{Gud-Sve-1}.

Let $W$ be an open subset of a non-compact semisimple Lie group
$G$ equipped with its standard Riemannian metric and
$\phi:W\to\cn$ be a real analytic map. Let $G^\cn$ denote the
complexification of the Lie group $G$. Then $\phi$ extends
uniquely to a holomorphic map $\phi^\cn:W^\cn\to\cn$ from some
open subset $W^\cn$ of $G^\cn$. By restricting this map to $U\cap
W^\cn$ we obtain a real analytic map $\phi^*:W^*\to\cn$ from some
open subset $W^*$ of $U$.

\begin{theorem}\label{theo:duality}
Let $\E$ be a family of maps $\phi:W\to\cn$ and $\E^*$ be the dual
family consisting of the maps $\phi^*:W^*\to\cn$ constructed as
above.  Then $\E^*$ is an eigenfamily if and only if $\E$ is an
eigenfamily.
\end{theorem}

\begin{proof}  Let $\g=\k+\p$ be a Cartan decomposition of the Lie
algebra of $G$ where $\k$ is the Lie algebra of a maximal compact
subgroup $K$.  Furthermore let the left-invariant vector fields
$X_1,\dots,X_n\in\p$ form a global orthonormal frame for the
distribution generated by $\p$ and similarly $Y_1,\dots,Y_m\in\k$
form a global orthonormal frame for the distribution generated by
$\k$. We shall now assume that $\phi$ and $\psi$ are elements of
the eigenfamily $\E$ on the Riemannian $W$ i.e.
$$\tau (\phi)=\sum_{k=1}^m Y_k^{2}(\phi) +\sum_{k=1}^n
X_k^{2}(\phi)=\lambda\cdot\phi,$$ $$\kappa
(\phi,\psi)=\sum_{k=1}^m Y_k(\phi)Y_k(\psi)+\sum_{k=1}^n
X_k(\phi)X_k(\psi)=\mu\cdot\phi\cdot\psi.$$ By construction and by
the unique continuation property of real analytic functions the
extension ${\phi}^\cn$ of $\phi$ satisfies the same equations.

The Lie algebra of $U$ has the decomposition $\un=\k+i\p$ and the
left-invariant vector fields $iX_1,\dots,iX_n\in\ i\p$ form a
global orthonormal frame for the distribution generated by $i\p$.
Then
$$\tau (\phi^*)=-\sum_{k=1}^m Y_k^{2}(\phi^*)
+\sum_{k=1}^n (iX_k)^{2}(\phi^*)=-\lambda\cdot\phi^*$$
$$\kappa (\phi^*,\phi^*)=-\sum_{k=1}^m Y_k(\phi^*)Y_k(\psi^*)
+\sum_{k=1}^n
(iX_k)(\phi^*)(iX_k)(\psi^*)=-\mu\cdot\phi^*\cdot\psi^*.$$ This
shows that $\E^*$ is an eigenfamily. The converse is similar.
\end{proof}

\begin{theorem}\label{theo:bi-duality}
Let $\B$ be a family of maps $\phi:W\to\cn$ and $\B^*$ be the dual
family consisting of the maps $\phi^*:W^*\to\cn$ constructed as
above.  Then $\B^*$ is a bi-eigenfamily if and only if $\B$ is a
bi-eigenfamily.
\end{theorem}

\begin{proof}
The argument is similar to that used for Theorem
\ref{theo:duality}
\end{proof}

\section{The semi-Riemannian Lie group $\GLC n$}

Let $h$ be the standard left-invariant semi-Riemannian metric on
the general linear group $\GLC n$ induced by the semi-Euclidean
scalar product on the Lie algebra $\glc n$ given by
$$h(Z,W)=\Re\trace ZW.$$ Then we have the orthogonal
decomposition $$\glc n=\W_+\oplus\W_-$$ of Lie algebra $\glc n$
where
$$\W_+=\{Z\in\glc n |\ Z-Z^*=0\}$$ is the subspace of Hermitian
matrices and
$$\W_-=\{Z\in\glc n |\ Z+Z^*=0\}$$ is the subspace of skew-Hermitian
matrices.  The scalar product is positive definite on $\W_+$ and
negative definite on $\W_-$.  This means that for two complex
valued functions $\phi,\psi:U\to\cn$ locally defined on $\GLC n$
the differential operators $\tau$ and the $\kappa$ satisfy
$$\tau(\phi)=\sum_{Z\in\B_+}Z^2(\phi)-\sum_{Z\in\B_-}Z^2(\phi),$$
$$\kappa(\phi,\psi)=\sum_{Z\in\B_+}Z(\phi)Z(\psi)
-\sum_{Z\in\B_-}Z(\phi)Z(\psi)$$ where $\B_+$ and $\B_-$ are
orthonormal bases for $\W_+$ and $\W_-$, respectively.

With the duality principle of Theorem \ref{theo:bi-duality} we can
now easily construct harmonic morphisms from the compact
Riemannian Lie groups
$$\SO n,\ \ \SU{n}\ \ \text{and}\ \ \Sp n$$
via the following classical dualities $G\cong U$:

$$\SLR n=\{x\in\GLR{n}\ |\ \det x=1\} \cong \SU n,$$

$$\SUs{2n}=\{g=\begin{pmatrix} z & w
\\ -\bar w & \bar z\end{pmatrix}
\ |\ g\in\SLC {2n}\}\cong \SU{2n},$$

$$\SOs{2n}=\{z\in\SUU nn\ |\ z\cdot I_{nn}\cdot J_{n}\cdot
z^t=I_{nn}\cdot J_{n}\}\cong\SO{2n},$$

$$\SpR n =\{g\in\SLR{2n}\ |\ g\cdot J_{n}\cdot
g^t=J_{n}\} \cong\Sp n.$$

$$\SOO pq=\{x\in\SLR{p+q}\ |\ x\cdot I_{pq}\cdot
x^t=I_{pq}\}\cong\SO{p+q},$$

$$\SUU pq=\{z\in\SLC{p+q}\ |\ z\cdot I_{pq}\cdot
z^*=I_{pq}\}\cong\SU{p+q},$$

$$\Spp pq=\{g\in\GLH{p+q}\ |\ g\cdot I_{pq}\cdot g^*
=I_{pq}\}\cong \Sp{p+q}.$$

\section{Acknowledgements}  The authors are grateful to Martin Svensson
for useful comment on this paper.

\appendix
\section{}\label{app:general}

In this section we prove the result stated in Theorem
\ref{theo:bi-rational}.  It shows how the elements of a
bi-eigenfamily $\B$ of complex valued functions on a
semi-Riemannian manifold $(M,g)$ can be used to produce a variety
of harmonic morphisms defined on open and dense subsets of $M$.
The first result shows how the operators $\tau$ and $\kappa$
behave with respect to products.

\begin{lemma}\label{lemm:products}
Let $(M,g)$ be a semi-Riemannian manifold and $\E_1,\E_2$ be two
families of complex valued functions on $M$.  If there exist
complex numbers
$\lambda_{1},\mu_{1},\lambda_{2},\mu_{2},\mu\in\cn$ such that for
all $\phi_1,\phi_2\in\E_1$ and $\psi_1,\psi_2\in\E_2$
$$\tau(\phi_1)=\lambda_{1}\phi_1,\ \ \kappa (\phi_1,\phi_2)
=\mu_{1} \phi_1\phi_2,$$ $$\tau(\psi_1)=\lambda_{2}\psi_1,\ \
\kappa (\psi_1,\psi_2)=\mu_{2} \psi_1\psi_2,$$ $$\kappa
(\phi_1,\psi_1)=\mu\phi_1\psi_1$$ then the following relations
hold
$$\tau (\phi_1\psi_1)=(\lambda_{1}+2\mu+\lambda_{2})\phi_1\psi_1,$$
$$\kappa (\phi_1\psi_1,\phi_2\psi_2)
=(\mu_{1}+2\mu+\mu_{2})\phi_1\psi_1\phi_2\psi_2$$
 for all $\phi_1,\phi_2\in\E_1$ and $\psi_1,\psi_2\in\E_2$.
\end{lemma}

\begin{proof}  The statement is an immediate consequence of the
following basic facts concerning first and second order
derivatives of products
$$X(\phi_1\psi_1)=X(\phi_1)\psi_1+\phi_1X(\psi_1),$$
$$X^2(\phi_1\psi_1)=X^2(\phi_1)\psi_1+2X(\phi_1)X(\psi_1)+\phi_1X^2(\psi_1).$$
\end{proof}

The following result shows how the operators $\tau$ and $\kappa$
behave with respect to quotients.

\begin{lemma}\label{lemm:quotients}
Let $(M,g)$ be a semi-Riemannian manifold and $P,Q:M\to\cn$ be two
complex valued functions on $M$. If there exists a complex number
$\lambda\in\cn$ such that $$\tau (P)=\lambda P\ \ \text{and}\ \
\tau (Q)=\lambda Q$$ then the quotient $\phi=P/Q$ is a harmonic
morphism if and only if
$$Q^2\kappa (P,P)=PQ\kappa (P,Q)=P^2\kappa (Q,Q).$$
\end{lemma}

\begin{proof}  For the first and second order derivatives of the quotient
$P/Q$ we have the following basic facts
$$X(\phi)=\frac{X(P)Q-PX(Q)}{Q^2},$$
$$X^2(\phi)=\frac{Q^2X^2(P)-2QX(P)X(Q)+2PX(Q)X(Q)-PQX^2(Q)}{Q^3}$$
leading to the following formulae for $\tau(\phi)$ and $\kappa
(\phi,\phi )$
$$Q^3\tau(\phi)=Q^2\tau(P)-2Q\kappa(P,Q)+2P\kappa(Q,Q)-PQ\tau(Q),$$
$$Q^4\kappa(\phi,\phi)=Q^2\kappa(P,P)-2PQ\kappa(P,Q)+P^2\kappa(Q,Q).$$
The statement is a direct consequence of those relations.
\end{proof}

\begin{proof}[Proof of Theorem \ref{theo:bi-rational}]
For the finite eigenfamilies $\E _1$ and $\E _2$ on $M$ we define
the infinite sequences
$$\{\E_1 ^k\}_{k=1}^\infty\ \ \text{and}\ \ \{\E_2 ^k\}_{k=1}^\infty$$
by induction

$$\E_1 ^1=\E_1,\ \
\E_1^{k+1}=\E_1^k\cdot\E_1=\{\phi\cdot\phi^1 |\ \phi\in\E_1^k,\
\phi^1\in\E_1\},$$

$$\E_2 ^1=\E_2,\ \
\E_2^{k+1}=\E_2^k\cdot\E_2=\{\psi\cdot\psi^1 |\ \psi\in\E_2^k,\
\psi^1\in\E_2\}.$$
It then follows from the fact that
$$\kappa (\phi,\phi^1 )=k\mu_1\phi\phi^1$$ for all
$\phi\in\E_1^k$ and $\phi^1\in\E_1^1$ and Lemma
\ref{lemm:products} that each $\E_1^{k+1}$ is an eigenfamily on
$M$.  The same is of course true for each $\E_2^{k+1}$.  A simple
calculation shows that for all $\Phi \in \E_1 ^k$, $\Psi \in \E_2
^l$
$$\kappa(\Phi, \Psi)=\mu\cdot k\cdot l \cdot\Phi \cdot\Psi.$$
With this at hand the statement of Theorem \ref{theo:bi-rational}
is an immediate consequence of Lemma \ref{lemm:quotients}.
\end{proof}


\begin{thebibliography}{99}

\bibitem{Bai-Eel}
P.~Baird and J.~Eells, {\it A conservation law for harmonic maps},
Geometry Symposium Utrecht 1980, Lecture Notes in Mathematics {\bf
894}, 1-25, Springer (1981).


\bibitem{Bai-Woo-book}
P.~Baird and J.~C. Wood, {\it Harmonic morphisms between
Riemannian manifolds}, London Math. Soc. Monogr. No. {\bf 29},
Oxford Univ. Press (2003).

\bibitem{Bai-Woo-1}
P.~Baird and J.~C. Wood, {\it Harmonic morphisms, Seifert fibre
spaces and conformal foliations}, Proc. London Math. Soc. {\bf 64}
(1992), 170-197.

\bibitem{Fug-1}
B.~Fuglede, {\it Harmonic morphisms between Riemannian manifolds},
Ann. Inst. Fourier {\bf 28} (1978), 107-144.

\bibitem{Fug-2} B.~Fuglede,
{\it Harmonic morphisms between semi-riemannian manifolds}, Ann.
Acad. Sci. Fennicae {\bf 21} (1996), 31-50.

\bibitem{Gud-bib}
S.~Gudmundsson, {\it The Bibliography of Harmonic
Morphisms}, {\tt http://www.matematik.lu.se/\\
matematiklu/personal/sigma/harmonic/bibliography.html}

\bibitem{Gud-1}
S.~Gudmundsson, {\it On the existence of harmonic morphisms from
symmetric spaces of rank one}, Manuscripta Math. {\bf 93} (1997),
421-433.

\bibitem{Gud-Sak-1}
S.~Gudmundsson and A.~Sakovich, {\it Harmonic morphisms from the
classical compact semisimple Lie groups}, preprint (Lund
University) 2006.

\bibitem{Gud-Sve-1}
S.~Gudmundsson and M.~Svensson, {\it Harmonic morphisms from the
Grassmannians and their non-compact duals}, Ann. Global Anal.
Geom. {\bf 30} (2006), 313-333.

\bibitem{Gud-Sve-2}
S.~Gudmundsson and M.~Svensson, {\it Harmonic morphisms from the
compact semisimple Lie groups and their non-compact duals},
Differ. Geom. Appl. 24 (2006), 351-366.

\bibitem{Gud-Sve-3}
S.~Gudmundsson and M.~Svensson, {\it On the existence of harmonic
morphisms from certain symmetric spaces}, J. Geom. Phys. {\bf 57}
(2007), 353-366.

\bibitem{Ish}
T.~Ishihara, {\it A mapping of Riemannian manifolds which
preserves harmonic functions}, J. Math. Kyoto Univ. {\bf 19}
(1979), 215-229.

\bibitem{Kna}
A.~W. Knapp, {\it Lie Groups Beyond an Introduction}, Birkh\"
auser (2002).




\end{thebibliography}
\end{document}